# Fractional Brownian fields, duality, and martingales

**Vladimir Dobrić[1] and Francisco M. Ojeda[2]**

*Lehigh University and Universidad Simón Bolívar and Lehigh University*

**Abstract:** In this paper the whole family of fractional Brownian motions is constructed as a single Gaussian field indexed by time and the Hurst index simultaneously. The field has a simple covariance structure and it is related to two generalizations of fractional Brownian motion known as multifractional Brownian motions. A mistake common to the existing literature regarding multifractional Brownian motions is pointed out and corrected. The Gaussian field, due to inherited "duality", reveals a new way of constructing martingales associated with the odd and even part of a fractional Brownian motion and therefore of the fractional Brownian motion. The existence of those martingales and their stochastic representations is the first step to the study of natural wavelet expansions associated to those processes in the spirit of our earlier work on a construction of natural wavelets associated to Gaussian-Markov processes.

## 1. Introduction

The basic quantities describing movements of a viscous fluid are related through the Reynolds number. When the Reynolds number exceeds a critical value, which depends on the fluid, the average velocity of the region and its geometry, the flow becomes unstable and random. In his work on the theory of turbulent flow Kolmogorov proposed a model [10] which assumes that the kinetic energy in large scale motions of a turbulent flow is transferred to smaller scale turbulent motions. At smaller scales the Reynolds number associated with that region is reduced. When the Reynolds number of a region falls below the critical value of the region, turbulent motion stops and the remaining kinetic energy is dissipated as heat. Kolmogorov assumed that smaller scale turbulent motions can be described by a random field. Kolmogorov modeling assumptions produce the random field $B_H$ that is self-similar with stationary increments and the second moment of its increments is of the form $\mathbb{E}\left|B_H(t+s) - B_H(s)\right|^2 = c\left|t\right|^{2H}$. Analyzing the properties of those fields [8, 9] Kolmogorov obtained a spectral representation for the fields with stationary increments.

In 1968 Mandelbrot and Van Ness interpreted the nonanticipating representation of $B_H$,

$$B_H(t) - B_H(s) = c \int_R \left((t-x)_+^{H-\frac{1}{2}} - (s-x)_+^{H-\frac{1}{2}}\right) db(x),$$

---

[1]Department of Mathematics, 14 Packer Av., Lehigh University, Bethlehem, PA 18015, e-mail: vd00@lehigh.edu

[2]Departamento de Matemáticas, Universidad Simón Bolívar, Apartado 89000, Caracas 1080-A, Venezuela, e-mail: fojeda@usb.ve

with respect to an orthogonal white noise $db$, obtained earlier by Pinsker and Yaglom, as a fractional integral, and called $B_H$ the fractional Brownian motion in the case when $b$ is a Brownian motion. Most of the above information is quoted from Molchan [14]. His paper is an excellent short review of the history of fractional Brownian motion before 1970. The index $H$ is sometimes called the Hurst exponent, after the British hydrologist H.E. Hurst, who studied the annual flows of the Nile.

Fractional Brownian motion has applications in financial mathematics, telecommunication networks, hydrology, physical and biological sciences, just to mention a few.

An appropriate representation of fractional Brownian motion is important for analyzing its properties and for computational purposes. Meyer, Sellan and Taqqu developed a method for constructing wavelet expansions of fractional Brownian motion [12]. Their construction was based on fractional integration and differentiation of wavelet expansions of a Gaussian white noise and encompasses a large class of wavelets.

We have found an iterative method for obtaining an orthogonal expansion of a Gaussian Markov process [5] directly from its covariance function. It turns out that our method produces a wavelet expansion if time is measured by the natural measure associated to the process. It is therefore natural to ask if it is possible to construct a "natural" wavelet expansion associated with a fractional Brownian motion (fBm in short) in the spirit of our work in [5]. The main properties used in our construction of wavelets associated with the Gaussian Markov processes were the invariance of the processes on pinning (stays Markov), independence of associated pinned processes, and the existence of associated martingales.

The first step toward finding natural wavelets for fBm is to investigate invariances for the entire class of fractional Brownian motions. That means considering the family of processes $(B_H)_{0<H<1}$ where $B_H = (B_H(t))_{t\in\mathbb{R}}$ is a fBm. Let us point out that it is not straight forward to prove that, for the fixed Hurst index $H$, the covariance function,

$$(1) \qquad \mathbb{E} B_H(t) B_H(s) = \frac{1}{2}(|t|^{2H} + |s|^{2H} - |t-s|^{2H}),$$

of fBm is indeed a positive definite form. It uses tricks and the fact that the characteristic function $\varphi$ of an $2H$ stable random variable is given by $\varphi(t) = e^{-a|t|^{2H}}$ for some constant $a$. In the search for invariances relevant to achieving our goal, that is to obtaining natural wavelet expansions for fBm, we have constructed a Gaussian field, indexed by $\mathbb{R} \times (0,1)$, that encompasses at once all fractional Brownian motions for all $H$. More precisely, this field has the property that, when the second parameter, $H$, is fixed, the resulting process is an $H$-fBm. We started our construction by using already mentioned nonanticipating representation of fBm, which, for $0 < H < 1$, can be written as

$$(2) \qquad B_H(t) = \frac{\sqrt{\Gamma(2H+1)\sin(\pi H)}}{\Gamma(H+\frac{1}{2})} \int_{-\infty}^{\infty} (t-x)_+^{H-\frac{1}{2}} - (-x)_+^{H-\frac{1}{2}}\, dW(x)$$

for $t \in R$, where $(W(t))_{t\in R}$ is a Brownian motion. When we started the computations our hope was that the covariance of our Gaussian field would resemble the covariance of fractional Brownian motion (equation (1)), with $2H$ substituted by $H + H'$ modulo a function of $H$ and $H'$. This does not turn out to be the case. The resulting covariance was more complex. Its form revealed that the nonanticipating representation may not be the most "efficient" representation of fBm. Section 3 of this paper contains all computations relevant to finding that covariance.

Although the Gaussian field developed in Section 3 did not meet our objective, its structure revealed what intrinsic dependences of the nonanticipating representation of fBm are, and what should be done to construct a new field having, for our purposes, the "right" covariance, namely the covariance

$$(3) \qquad \mathbb{E} B_H(t) B_{H'}(s) = \alpha_{H,H'}(|t|^{H+H'} + |s|^{H+H'} - |t-s|^{H+H'}).$$

We have succeeded to obtain such a field. Furthermore we have found a nonanticipating representation of $\{B_H(t)\}_{t \in \mathbb{R}, H \in (0,1)}$, which, when $H$ is fixed, equals in distribution to the standard nonanticipating representation of fBm. Those results are presented in Section 4.

Our new field has the property that when $H + H' = 1$ the right hand side of (3) becomes the covariance of a Brownian motion. In that case we call $B_H, B_{H'}$ a dual pair. A particular property of a dual pair is that it generates two martingales, one driving $B_H$ and the other $B_{H'}$. In [4] we have obtained stochastic integral representations for those martingales. It turns out that our representation coincides with the fundamental martingale of a fBm discovered by Molchan and Golosov [13, 15]. Our work on that subject is in Section 5.

Lévy Véhel and Peltier [23] and Benassi et al. [2] have independently introduced a multifractional Brownian motion. In Section 6 we clarify the relationship between multifractional Brownian motion and the fractional Brownian fields introduced in Section 3 and 4. Furthermore we point out and correct an error in the covariance of the multifractional Brownian motion obtained from the nonanticipating moving average of fBm and show that in fact the processes of Lévy Véhel and Peltier [23] and Benassi et al. [2] are not the same, as it has been claimed in Cohen [3].

For readers convenience some basic facts about fractional Brownian motions and the notation that is used in the paper is included in Section 2.

## 2. Notation and preliminaries

This section is a brief overview of some properties of the fractional Brownian motion. The material presented here can be found in Taqqu [22], Molchan [14], chapter 7 in Samorodnitsky and Taqqu [21] and Dzhaparidze and Van Zanten [6]. All the processes and fields considered will be real valued.

The following proposition summarizes some properties of $H$-self-similar processes with stationary increments and finite second moments.

**Proposition 1.** *If $Z = \{Z_t\}_{t \in \mathbb{R}}$ is $H$-self-similar with stationary increments and finite second moments then*

1. $Z_0 = 0$ *with probability one.*
2. *If $H \neq 1$, then $\mathbb{E} Z_t = 0$ for all $t \in \mathbb{R}$.*
3. $H \leq 1$.
4. *The covariance of $Z$ is given by*

$$(4) \qquad \mathbb{E}(Z_s Z_t) = \frac{\mathbb{E}(Z_1^2)}{2}\left(|s|^{2H} + |t|^{2H} - |s-t|^{2H}\right).$$

As already mentioned in the introduction Kolmogorov described an $H$-self-similar process with stationary increments and its covariance (4) in [9]. A centered Gaussian process $\{Z_t\}_{t \in \mathbb{R}}$ with covariance (4) is called a fractional Brownian motion (fBm) with index $H$. If $\mathbb{E}(Z_1^2) = 1$, the process is a standard fBm. Since the distribution of a Gaussian process is entirely determined by its covariance, (4) implies that fBm is $H$-self-similar and has stationary increments. Conversely, an $H$-self-similar Gaussian process with stationary increments is a fBm. For $H = \frac{1}{2}$ a standard fractional Brownian motion is a standard Brownian motion on the real line, that is,

$$B_{\frac{1}{2}}(t) = \begin{cases} B^1(-t) & t < 0 \\ B^2(t) & t \geq 0 \end{cases},$$

where $\{B_t^1\}_{t \geq 0}$, $\{B_t^2\}_{t \geq 0}$ are two independent standard Brownian motions.

Unless otherwise stated all Brownian motions and fractional Brownian motions appearing on this work will be assumed to be standard.

Using Kolmogorov-Centsov criterion for Hölder continuity of a process, for every $\gamma$ in $(0, H)$, an $H$-fBm has a continuous modification that is locally Hölder continuous with exponent $\gamma$. Throughout the paper we will assume that fBm is such a modification.

Let $\{B_t\}_{t \in R}$ and $\{W_t\}_{t \in R}$ be Brownian motions on the real line. Then the process defined by the moving average stochastic integral

$$(5) \qquad B_H(t) = \frac{1}{c_H} \int_{-\infty}^{\infty} (t-x)_+^{H-\frac{1}{2}} - (-x)_+^{H-\frac{1}{2}} \, dB_x \text{ for } t \in \mathbb{R}$$

is a fBm, where

$$(6) \qquad c_H = \frac{\Gamma\left(H + \frac{1}{2}\right)}{\sqrt{\Gamma(2H+1)\sin(\pi H)}}$$

and $x_+^\alpha = 0$ if $x \leq 0$ and $x_+^\alpha = x^\alpha$ if $x > 0$. This representation is called nonanticipating, since for $t \geq 0$, it involves only integration over $(-\infty, t]$. Another moving average representation of fBm for $H \neq \frac{1}{2}$ is given by

$$(7) \qquad W_H(t) = \frac{1}{d_H} \int_{-\infty}^{\infty} |t-x|^{H-\frac{1}{2}} - |x|^{H-\frac{1}{2}} \, dW_x \text{ for } t \in \mathbb{R}$$

where

$$(8) \qquad d_H = \frac{\Gamma\left(H + \frac{1}{2}\right)}{\sqrt{\Gamma(2H+1)}} \sqrt{2 \left(\frac{1 - \sin(\pi H)}{\sin(\pi H)}\right)},$$

and when $H = \frac{1}{2}$ by

$$(9) \qquad W_{\frac{1}{2}}(t) = \frac{1}{d_{\frac{1}{2}}} \int_{-\infty}^{\infty} \log\left(\frac{1}{|t-x|}\right) - \log\left(\frac{1}{|x|}\right) dW_x \text{ for } t \in \mathbb{R}$$

where $d_{\frac{1}{2}} = \pi$. This representation is called well-balanced.

The odd and even part of fBm play an important role in this paper. If $\{Z_H(t)\}_{t \in \mathbb{R}}$ is a $H$-fBm, then the odd and the even process of fBm are defined by

$$Z_H^o(t) = \frac{1}{2}(Z_H(t) - Z_H(-t)), \text{ and } Z_H^e(t) = \frac{1}{2}(Z_H(t) + Z_H(-t)), \; t \in \mathbb{R}$$

respectively. It is straightforward to check that the processes $Z_H^o$ and $Z_H^e$ are independent and that their covariances are given by

$$(10) \qquad \mathbb{E} Z_H^o(s) Z_H^o(t) = \frac{1}{4}\left(|s+t|^{2H} - |s-t|^{2H}\right)$$

and

$$\mathbb{E} Z_H^e(s) Z_H^e(t) = \frac{|s|^{2H} + |t|^{2H}}{2} - \frac{|s+t|^{2H} + |s-t|^{2H}}{4} \tag{11}$$

for $s, t$ in $\mathbb{R}$. Both $Z_H^e$ and $Z_H^o$ are $H$-self-similar processes and by assumption have continuous sample paths with probability 1. Clearly, each path of the process $Z_H^o$ ($Z_H^e$) is an odd (even) function on a set of probability 1, and so it suffices to consider these processes for $t \geq 0$. For $H = \frac{1}{2}$ the odd and even part of fBm are Brownian motions (up to a constant multiple).

Moving average representations for the odd and even part of fBm can be found in Dzhaparidze and Van Zanten [6], Nuzman and Poor [17].

## 3. Dependent fractional Brownian field

In this section we will assume that each element of the family $\{B_H(t)\}_{t \in \mathbb{R}, H \in (0,1)}$ is represented by a nonanticipating moving average stochastic integral (5) respect to the same $\{B_t\}_{t \in \mathbb{R}}$. For $0 < H, H' < 1$ set

$$c_{H,H'} = \frac{\sqrt{\Gamma(2H+1)\sin(\pi H)}\sqrt{\Gamma(2H'+1)\sin(\pi H')}}{\pi}.$$

**Theorem 2.** *Let $K$ be the covariance $\mathbb{E} B_H(t) B_{H'}(s)$. If $H + H' \neq 1$ then*

$$\begin{aligned}
K = c_{H,H'} \Gamma\left(-\left(H + H'\right)\right) &\left\{ \cos\left(\left(H' - H\right)\frac{\pi}{2}\right) \cos\left(\left(H + H'\right)\frac{\pi}{2}\right) \right. \\
&\times \left(|t-s|^{H+H'} - |t|^{H+H'} - |s|^{H+H'}\right) \\
&- \sin\left(\left(H' - H\right)\frac{\pi}{2}\right) \sin\left(\left(H + H'\right)\frac{\pi}{2}\right) \\
&\left. \times \left(sgn(t)|t|^{H+H'} - sgn(s)|s|^{H+H'} - sgn(t-s)|t-s|^{H+H'}\right) \right\}.
\end{aligned} \tag{12}$$

*If $H + H' = 1$ and $t \neq s$, $t \neq 0$, $s \neq 0$, then*

$$\begin{aligned}
K = c_{H,H'} &\left\{ \cos\left(\left(H' - H\right)\frac{\pi}{2}\right) \frac{\pi}{2} (|t| + |s| - |t-s|) \right. \\
&- \sin\left(\left(H' - H\right)\frac{\pi}{2}\right) \\
&\left. \times (t \log|t| - s \log|s| - (t-s) \log|t-s|) \right\}.
\end{aligned} \tag{13}$$

*If $H + H' = 1$ and $t = s$, $t \neq 0$, then*

$$K = c_{H,H'} \cos\left(\left(H' - H\right)\frac{\pi}{2}\right) \frac{\pi}{2} \{|t| + |s|\}. \tag{14}$$

*If $H + H' = 1$ and $t = 0$ or $s = 0$, then*

$$K = 0. \tag{15}$$

Fourier transform is the main tool in proving the theorem. Before proving it we will establish some technical results.

Let $\widehat{f}(\xi) = \int_{-\infty}^{\infty} e^{i\xi x} f(x) \, dx$ be the Fourier transform of a function that belong to $L^1(\mathbb{R}) \cap L^2(\mathbb{R})$ With this convention $\sqrt{2\pi} \|f\|_2 = \left\|\widehat{f}\right\|_2$ for all $f \in L^2(\mathbb{R})$.

**Lemma 3.** *For $0 < H < 1$ and $t \in \mathbb{R}$*

$$\text{(16)} \quad \left(\widehat{\frac{(t-\cdot)_+^{H-\frac{1}{2}} - (-\cdot)_+^{H-\frac{1}{2}}}{\Gamma\left(H+\frac{1}{2}\right)}}\right)(\xi) = \frac{e^{it\xi} - 1}{i\xi} (i\xi)^{-\left(H-\frac{1}{2}\right)},$$

$$\text{(17)} \quad \left(\widehat{\frac{(t-\cdot)_-^{H-\frac{1}{2}} - (-\cdot)_-^{H-\frac{1}{2}}}{\Gamma\left(H+\frac{1}{2}\right)}}\right)(\xi) = -\frac{e^{it\xi} - 1}{i\xi} (-i\xi)^{-\left(H-\frac{1}{2}\right)},$$

*where*

$$\text{(18)} \quad (i\xi)^{-\left(H-\frac{1}{2}\right)} = \begin{cases} |\xi|^{-\left(H-\frac{1}{2}\right)} e^{i\frac{\pi}{2}\left(H-\frac{1}{2}\right)} & \xi < 0 \\ |\xi|^{-\left(H-\frac{1}{2}\right)} e^{-i\frac{\pi}{2}\left(H-\frac{1}{2}\right)} & \xi > 0 \end{cases}.$$

*Proof.* Equation (16) follows from the theory of fractional integration and differentiation, see for example Lemma 1 and Lemma 3 in [20]. Let

$$f_t(x) = \frac{(t-x)_+^{H-\frac{1}{2}} - (-x)_+^{H-\frac{1}{2}}}{\Gamma\left(H+\frac{1}{2}\right)},$$

then the identity

$$\frac{(t-x)_-^{H-\frac{1}{2}} - (-x)_-^{H-\frac{1}{2}}}{\Gamma\left(H+\frac{1}{2}\right)} = f_{-t}(-x)$$

is a consequence of the relation $x_+^\alpha = (-x)_-^\alpha$. Equation (17) now follows from (16). □

**Lemma 4.** *For $s, t \in \mathbb{R}$ and $0 < H, H' < 1$ let*

$$\text{(19)} \quad I_1 \stackrel{def}{=} \int_0^\infty \frac{\sin^2\left(\frac{t\xi}{2}\right) + \sin^2\left(\frac{s\xi}{2}\right) - \sin^2\left(\frac{\xi(t-s)}{2}\right)}{\xi^{1+H+H'}} d\xi.$$

*If $H + H' \neq 1$ then*

$$\text{(20)} \quad \begin{aligned} I_1 &= \frac{\Gamma\left(-\left(H+H'\right)\right)\cos\left(\left(H+H'\right)\frac{\pi}{2}\right)}{2} \\ &\quad \times \left(|t-s|^{H+H'} - |t|^{H+H'} - |s|^{H+H'}\right) \end{aligned}$$

*and if $H + H' = 1$ then*

$$I_1 = \frac{\pi}{4}\left(|t| + |s| - |t-s|\right).$$

*Proof.* The proof is trivial if either $s = 0$ or $t = 0$, so assume that $s \neq 0$, $t \neq 0$. Formula 3.823 in Gradshteyn and Ryzhik [7] states

$$\text{(21)} \quad \int_0^\infty x^{\mu-1} \sin^2(ax)\, dx = -\frac{\Gamma(\mu)\cos\left(\frac{\mu\pi}{2}\right)}{2^{\mu+1} a^\mu} \text{ for } a > 0 \text{ and } -2 < \text{Re}(\mu) < 0.$$

When $H + H' \neq 1$ set $\mu = -\left(H + H'\right)$. Observe that $\sin^2(x) = \sin^2(|x|)$. Applying the identity (21) to the right hand side of (19) when $s = t$ gives

$$I_1 = -\Gamma\left(-H - H'\right)\cos\left(\left(H + H'\right)\frac{\pi}{2}\right)|t|^{H+H'},$$

and when $s \neq t$

$$I_1 = \frac{\Gamma\left(-H - H'\right)}{2}\cos\left(\left(H + H'\right)\frac{\pi}{2}\right)\beta(H + H')$$
$$\times \left\{|t - s|^{H+H'} - |t|^{H+H'} - |s|^{H+H'}\right\}.$$

In the case of $H + H' = 1$ we will use formula 3.821 (9) from Gradshteyn and Ryzhik [7],

$$\int_0^\infty \frac{\sin^2(ax)}{x^2}dx = \frac{a\pi}{2} \text{ for } a > 0.$$

to conclude that when $s \neq t$ then $I_1 = \pi/4\left(|t| + |s| - |t - s|\right)$, and when $s = t$ then $I_1 = 1/2\,|t|\,\pi$. □

**Lemma 5.** *For $s, t \in \mathbb{R}$ and $0 < H, H' < 1$ let*

(22) $$I_2 \stackrel{def}{=} \int_0^\infty \frac{\sin(\xi(t - s)) + \sin(s\xi) - \sin(t\xi)}{\xi^{1+H+H'}}d\xi.$$

*If $s = t$ or $s = 0$ or $t = 0$ then $I_2 = 0$. Otherwise, if $H + H' \neq 1$ then*

(23) $$I_2 = \Gamma\left(-\left(H + H'\right)\right)\sin\left(\left(H + H'\right)\frac{\pi}{2}\right)$$
$$\times \left\{sgn(t)|t|^{H+H'} - sgn(s)|s|^{H+H'} - sgn(t - s)|t - s|^{H+H'}\right\},$$

*and if $H + H' = 1$ then*

(24) $$I_2 = t\log|t| - s\log|s| - (t - s)\log|t - s|.$$

*Proof.* The proof is trivial if either $s = 0$ or $t = 0$ or $s = t$, so assume that $s \neq 0$, $t \neq 0$ and $s \neq t$. Formula 3.761 (4) in Gradshteyn and Ryzhik [7] reads:

(25) $$\int_0^\infty x^{\mu-1}\sin(ax)\,dx = \frac{\Gamma(\mu)}{a^\mu}\sin\left(\frac{\mu\pi}{2}\right) \text{ for } a > 0 \text{ and } 0 < |\text{Re}(\mu)| < 1.$$

Set $\mu = -\left(H + H'\right)$, and observe that $\sin(\xi x) = sgn(x)\sin(\xi|x|)$ Applying (25) to the right hand side of (22) in the case when $0 < H + H' < 1$ yields (23). In the case when $1 < H + H' < 2$, by the dominated convergence theorem, it follows that

$$I_2 = \lim_{\substack{a \to 0^+ \\ b \to \infty}} \int_a^b \frac{sgn(t-s)\sin(\xi|t-s|) + sgn(s)\sin(|s|\xi) - sgn(t)\sin(|t|\xi)}{\xi^{1+H+H'}} d\xi$$

$$= \lim_{\substack{a \to 0^+ \\ b \to \infty}} \frac{1}{H+H'} \int_a^b (sgn(t-s)|t-s|\cos(\xi|t-s|)$$

$$+ sgn(s)|s|\cos(|s|\xi) - sgn(t)|t|\cos(|t|\xi))\xi^{-H-H'} d\xi$$

$$= \frac{1}{(H+H')(H+H'-1)} \int_0^\infty \Big(-sgn(t-s)|t-s|^2 \sin(\xi|t-s|)$$

$$- sgn(s)|s|^2 \sin(|s|\xi) + sgn(t)|t|^2 \sin(|t|\xi)\Big)\xi^{1-H-H'} d\xi,$$

where the last two equalities are result of integration by parts and the fact that the boundary terms converge to zero as $b \to \infty$, and $a \to 0^+$. Applying (25) with $\mu = 2 - (H + H')$ and using $x\Gamma(x) = \Gamma(x+1)$ the equation (23) now follows readily.

Let us turn to the case $H + H' = 1$. For $x \in [1, 1.5]$ set

$$f(x) = \int_0^\infty \frac{sgn(t-s)\sin(\xi|t-s|) + sgn(s)\sin(|s|\xi) - sgn(t)\sin(|t|\xi)}{\xi^{1+x}} d\xi.$$

When $x = H+H'$ the function $f$ equals to the right hand side of (23), The integrand that defines $f$ is bounded by

$$g(\xi) = \begin{cases} \frac{|sgn(t-s)\sin(\xi|t-s|) + sgn(s)\sin(|s|\xi) - sgn(t)\sin(|t|\xi)|}{\xi^{1+1.5}} & \xi \in (0,1) \\ \frac{3}{\xi^{1+1}} & \xi \in [1,\infty) \end{cases}$$

which is an integrable function. By the dominated convergence theorem $f$ is continuous on $[1, 1.5]$. Finally $\lim_{x \downarrow 1} f(x)$ establishes (24), where we used (23) and the gamma function property $y\Gamma(y) = \Gamma(y+1)$ to rewrite $f(x)$ for $x \in (1, 1.5]$ as

$$f(x) = \frac{\Gamma(2-x)}{(-x)(1-x)} \sin\left(\frac{x\pi}{2}\right) \{sgn(t)|t|^x - sgn(s)|s|^x - sgn(t-s)|t-s|^x\}$$

and L'Hospital rule to compute the limit

$$\lim_{x \downarrow 1} \frac{sgn(t)|t|^x - sgn(s)|s|^x - sgn(t-s)|t-s|^x}{(1-x)}.$$

□

We have prepared the groundwork to prove Theorem 2.

*Proof.* By the Ito isometry

$$K = \frac{1}{c_H c_{H'}} \int_{-\infty}^\infty \left((t-x)_+^{H-\frac{1}{2}} - (-x)_+^{H-\frac{1}{2}}\right)\left((s-x)_+^{H'-\frac{1}{2}} - (-x)_+^{H'-\frac{1}{2}}\right) dx,$$

and by Plancherel identity and Lemma 3

$$K = \frac{\Gamma\left(H + \frac{1}{2}\right)\Gamma\left(H' + \frac{1}{2}\right)}{2\pi c_H c_{H'}} \int_{-\infty}^{\infty} \frac{e^{it\xi} - 1}{i\xi}(i\xi)^{-(H-\frac{1}{2})} \overline{\frac{e^{is\xi} - 1}{i\xi}(i\xi)^{-(H'-\frac{1}{2})}} d\xi$$

$$= \frac{\Gamma\left(H + \frac{1}{2}\right)\Gamma\left(H' + \frac{1}{2}\right)}{2\pi c_H c_{H'}} \left\{ \int_{-\infty}^{0} \frac{(e^{it\xi} - 1)(e^{-is\xi} - 1)}{|\xi|^{1+H+H'}} e^{-i(H'-H)\frac{\pi}{2}} d\xi \right.$$

$$\left. + \int_{0}^{\infty} \frac{(e^{it\xi} - 1)(e^{-is\xi} - 1)}{|\xi|^{1+H+H'}} e^{i(H'-H)\frac{\pi}{2}} d\xi \right\}.$$

Substituting $\xi' = -\xi$ in the first integral of the last equality above and then combining the two integrals yields

$$(26) \quad K = \frac{\Gamma\left(H + \frac{1}{2}\right)\Gamma\left(H' + \frac{1}{2}\right)}{\pi c_H c_{H'}} \operatorname{Re}\left(e^{i(H'-H)\frac{\pi}{2}} \int_{0}^{\infty} \frac{(e^{it\xi} - 1)(e^{-is\xi} - 1)}{|\xi|^{1+H+H'}} d\xi\right).$$

Using Euler's formula on $e^{i(H'-H)\frac{\pi}{2}}$ and $(e^{it\xi} - 1)(e^{-is\xi} - 1)$ and the identity $\sin^2 x = \frac{1 - \cos 2x}{2}$ we obtain

$$\operatorname{Re}\left(e^{i(H'-H)\frac{\pi}{2}} \int_{0}^{\infty} \frac{(e^{it\xi} - 1)(e^{-is\xi} - 1)}{|\xi|^{1+H+H'}} d\xi\right)$$

$$(27) \quad = \cos\left((H' - H)\frac{\pi}{2}\right) 2 \int_{0}^{\infty} \frac{\sin^2\left(\frac{t\xi}{2}\right) + \sin^2\left(\frac{s\xi}{2}\right) - \sin^2\left(\frac{\xi(t-s)}{2}\right)}{\xi^{1+H+H'}} d\xi$$

$$- \sin\left((H' - H)\frac{\pi}{2}\right) \int_{0}^{\infty} \frac{\sin(\xi(t-s)) + \sin(s\xi) - \sin(t\xi)}{\xi^{1+H+H'}} d\xi.$$

Observing that

$$(28) \quad \frac{\Gamma\left(H + \frac{1}{2}\right)\Gamma\left(H' + \frac{1}{2}\right)}{c_H c_{H'} \pi} = \frac{\sqrt{\Gamma(2H+1)\sin(\pi H)}\sqrt{\Gamma(2H'+1)\sin(\pi H')}}{\pi},$$

the expressions for $\mathbb{E} B_H(t) B_{H'}(s)$ now follow from equations (26), (27), (28) and Lemmas 4, 5 □

We will call a centered Gaussian field $\{B_H(t)\}_{t\in\mathbb{R}, H\in(0,1)}$ with the covariance given by Theorem 2 a dependent fractional Brownian field and refer to it as dfBf. The rest of the section elaborates on a property of the field that justifies that name.

Let $\{B_H^o(t)\}_{t\in[0,\infty), H\in(0,1)}$ and $\{B_H^e(t)\}_{t\in[0,\infty), H\in(0,1)}$, be the odd and even part of the dfBf $\{B_H(t)\}_{t\in\mathbb{R}, H\in(0,1)}$, that is

$$B_H^o(t) = \frac{B_H(t) - B_H(-t)}{2} \text{ and } B_H^e(t) = \frac{B_H(t) + B_H(-t)}{2}, t \geq 0.$$

For $H + H' \neq 1$ set

$$(29) \quad a_{H,H'} = -2\frac{\sqrt{\Gamma(2H+1)\sin(\pi H)}\sqrt{\Gamma(2H'+1)\sin(\pi H')}}{\pi}$$
$$\times \Gamma\left(-(H+H')\right)\cos\left((H'-H)\frac{\pi}{2}\right)\cos\left((H+H')\frac{\pi}{2}\right)$$

and for $H + H' = 1$

$$a_{H,H'} = \sqrt{\Gamma(2H+1)\Gamma(3-2H)} \sin^2(\pi H). \tag{30}$$

**Theorem 6.** *The covariance of the odd part and the even part of dfBf are*

$$\mathbb{E} B_H^o(t) B_{H'}^o(s) = a_{H,H'} \frac{|t+s|^{H+H'} - |t-s|^{H+H'}}{4} \tag{31}$$

*and*

$$\begin{aligned} &\mathbb{E} B_H^e(t) B_{H'}^e(s) \\ &= a_{H,H'} \left( \frac{|t|^{H+H'} + |s|^{H+H'}}{2} - \frac{|t-s|^{H+H'} + |t+s|^{H+H'}}{4} \right) \end{aligned} \tag{32}$$

*respectively.*

*Proof.* The covariance of the dfBf (Theorem 2), is of the form

$$\begin{aligned} \mathbb{E} B_H(t) B_{H'}(s) &= f(H, H') \left( |t|^{H+H'} + |s|^{H+H'} - |t-s|^{H+H'} \right) \\ &\quad + g(s, t, H, H'). \end{aligned} \tag{33}$$

It is a matter of straightforward computation to check that in both cases

$$\mathbb{E} B_H(t) B_{H'}(s) + \mathbb{E} B_H(-t) B_{H'}(-s) \text{ and } \mathbb{E} B_H(t) B_{H'}(-s) + \mathbb{E} B_H(-t) B_{H'}(s)$$

the $g$ function cancels. The result of the theorem follows by simple algebraic manipulation of the first part of the right hand side of (33) only. $\square$

The Gaussian fields $\{B_H^o(t)\}_{t \in [0,\infty), H \in (0,1)}$ and $\{B_H^e(t)\}_{t \in [0,\infty), H \in (0,1)}$ with covariances given by Theorem 6, will be called the odd and the even fractional Brownian field respectively. It is very simple to check that for every $a > 0$

$$\{B_H^o(at)\}_{t \in [0,\infty), H \in (0,1)} \stackrel{f.d.d.}{=} \{a^H B_H^o(at)\}_{t \in [0,\infty), H \in (0,1)}$$

and

$$\{B_H^e(at)\}_{t \in [0,\infty), H \in (0,1)} \stackrel{f.d.d.}{=} \{a^H B_H^e(at)\}_{t \in [0,\infty), H \in (0,1)},$$

where $\stackrel{f.d.d.}{=}$ indicates the equality of finite dimensional distributions.

Given a fBm its odd and even part are independent processes (indexed by $t$). However, this is not the case with the dfBf $\{B_H(t)\}_{t \in \mathbb{R}, H \in (0,1)}$. The fields $\{B_H^o(t)\}_{t \in [0,\infty), H \in (0,1)}$ and $\{B_H^e(t)\}_{t \in [0,\infty), H \in (0,1)}$ are not independent. For example if $H + H' = 1$ then

$$\begin{aligned} \mathbb{E} B_H^e(t) B_{H'}^o(s) &= \frac{1}{\pi} \sqrt{\Gamma(2H+1) \sin(\pi H)} \\ &\quad \times \sqrt{\Gamma(2H'+1) \sin(\pi H')} \sin\left((H'-H)\frac{\pi}{2}\right) \\ &\quad \times \left\{ s \log |s| - \frac{(t+s) \log|t+s| - (t-s) \log|t-s|}{2} \right\}. \end{aligned}$$

which is clearly not equal to 0. That is the reason for calling that field the dependent fractional Brownian field.

Another glance at the computation of the covariance $\mathbb{E}B_H^i(t)B_{H'}^j(s)$, $i,j \in \{o,e\}$, reveals that when $i = j$ the $g$ part of (33) cancels out while in the case when $i \neq j$ the first part of (33) cancels out leaving the $g$ part. Therefore the existence the $g$ part in (33) is the reason for dependence between $\{B_H^o(t)\}_{t\in\mathbb{R},H\in(0,1)}$ and $\{B_H^e(t)\}_{t\in\mathbb{R},H\in(0,1)}$. So it is natural to search for a method of creating a fractional Brownian field that would be of the form (33) with $g = 0$. One way of attacking that problem is a direct verification of positive definiteness of such a form, a very unattractive task. In the next section we present a straightforward construction of the field with the desired covariance.

## 4. Fractional Brownian field

The last remark of the previous points the direction for constructing a fractional Brownian field $\{B_H(t)\}_{t\in\mathbb{R},H\in(0,1)}$ with the covariance of the form of (33) with $g = 0$. The new Gaussian field contains all fractional Brownian motions too.

**Theorem 7.** *Let $B = \{B_H(t)\}_{t\in\mathbb{R},H\in(0,1)}$ and $W = \{W_H(t)\}_{t\in\mathbb{R},H\in(0,1)}$ be two dfBf generated by two independent Brownian motions $\{B_t\}_{t\in\mathbb{R}}$ and $\{W_t\}_{t\in\mathbb{R}}$ respectively. Let $\{B_H^i(t)\}_{t\in[0,\infty),H\in(0,1)}$, $i = o$ be the odd and $i = e$ be the even part of $B$, and let $\{W_H^i(t)\}_{t\in[0,\infty),H\in(0,1)}$, $i = o$ be the odd and $i = e$ the even part of $W$. Then the fractional Brownian field $\{Z_H(t)\}_{t\in\mathbb{R},H\in(0,1)}$ defined by*

$$(34) \qquad Z_H(t) = \begin{cases} B_H^e(t) + W_H^o(t) & \text{for } t \geq 0 \\ B_H^e(-t) - W_H^o(-t) & \text{for } t < 0 \end{cases}$$

*has the covariance*

$$(35) \qquad \mathbb{E}Z_H(t)Z_{H'}(s) = a_{H,H'}\left\{\frac{|t|^{H+H'} + |s|^{H+H'} - |t-s|^{H+H'}}{2}\right\},$$

*where $a_{H,H'}$ is given by equations (29) and (30).*

*Proof.* The proof follows from (32), (31) and independence of $\{B_t\}_{t\in\mathbb{R}}$ and $\{W_t\}_{t\in\mathbb{R}}$. □

We will call the process $\{Z_H(t)\}_{t\in\mathbb{R},H\in(0,1)}$ fractional Brownian field (fBf in short). Note that for any $t \in \mathbb{R}$,

$$Z_H(t) = \frac{1}{2}(B_H(t) + W_H(t) + B_H(-t) - W_H(-t)),$$

and that $\left\{\frac{B_H(t)+W_H(t)}{\sqrt{2}}\right\}_{t\in\mathbb{R}}$ and $\left\{\frac{B_H(t)-W_H(t)}{\sqrt{2}}\right\}_{t\in\mathbb{R},H\in(0,1)}$ are two independent dfBf. Consequently

$$(36) \qquad Z_H(t) = \frac{X_H(t) + Y_H(-t)}{\sqrt{2}}$$

where $\{X_H(t)\}_{t\in\mathbb{R},H\in(0,1)}$ and $\{Y_H(t)\}_{t\in\mathbb{R},H\in(0,1)}$ are two independent fractional Brownian fields, that is $\{Z_H(t)\}_{t\in\mathbb{R},H\in(0,1)}$ is a properly symmetrized dfBf.

**Proposition 8.** *Let $\{B_t\}_{t\in\mathbb{R}}$ and $\{W_t\}_{t\in\mathbb{R}}$ be two independent Brownian motion processes on the real line and let $\{Z_H(t)\}_{t\in\mathbb{R},H\in(0,1)}$ be defined by (34). Then*

$$Z_H(t) = \frac{1}{\sqrt{2}c_H}\int_{-\infty}^{\infty} (t-x)_+^{H-\frac{1}{2}} - (-x)_+^{H-\frac{1}{2}}\, dB_x$$
$$+ \frac{1}{\sqrt{2}c_H}\int_{-\infty}^{\infty} (-t-x)_+^{H-\frac{1}{2}} - (-x)_+^{H-\frac{1}{2}}\, dW_x.$$

*Proof.* Follows directly from (36) and (5). □

A fBf has the same self-similarity property in the time variable as the odd and even fractional Browian field, namely for $a > 0$

$$\{Z_H(at)\}_{t\in\mathbb{R},H\in(0,1)} \stackrel{f.d.d.}{=} \{a^H Z_H(at)\}_{t\in\mathbb{R},H\in(0,1)}.$$

Moreover, the stationary in the time variable of increments of the fBf easily follows from Theorem 7, that is

$$\{W_H(t)\}_{t\in\mathbb{R},H\in(0,1)} \stackrel{f.d.d.}{=} \{W_H(t+\delta) - W_H(\delta)\}_{t\in\mathbb{R},H\in(0,1)}$$

for any $\delta$.

An immediate consequence of the covariance structure of a fractional Brownian field is that when $H + H' = 1$ then

(37) $$\mathbb{E}(Z_H(t) Z_{H'}(t)) = a_{H,H'} \begin{cases} |s|\wedge|t| & \text{for } 0 \leq s,t \text{ or } 0 \geq s,t \\ 0 & \text{otherwise} \end{cases}.$$

That property leads to a construction of martingales associated to fractional Brownian motions. The methodology of the construction is the subject of the next section.

## 5. Duality and fundamental martingales

In what follows it is assumed that $\{B_H(t)\}_{t\in\mathbb{R},H\in(0,1)}$ is an fBf. Whenever $H+H' = 1$ we will call $B_H, B_{H'}$ (or $B_H^o, B_{H'}^o$ or $B_H^e, B_{H'}^e$) a dual pair. Dual pairs have unique properties. They generate martingales associated in a natural way to fractional Brownian motions $B_H$ and $B_{H'}$. The construction and explanation of the nature of those martingales is the subject of this section.

Every fBf is a sum of an even and an odd part of two independent dfBf's. For that reason it suffices to construct martingales, $M_H^o$ and $M_H^e$, adapted to the filtrations of the odd and even part of $\{B_H(t)\}_{t\in\mathbb{R}}$ respectively. The filtration generated by $M_H^o$ ($M_H^e$) coincides with the filtration of the odd (even) part of fBm, and for that reason, following the terminology used in Norros et al. [16] to describe a martingale for the fBm originally discovered by Molchan and Golosov (see Molchan [14]), we call $M_H^o$ ($M_H^e$) a fundamental martingale for the odd (even) part of fBm. Furthermore we derive a stochastic integral representation for those martingales. In a similar fashion this was done in Pipiras and Taqqu [18] and Pipiras and Taqqu [19] for the fractional Brownian motion.

For $i \in \{o, e\}$ set

$$\mathcal{F}_t^{H,i} = \sigma\left(B_H^i(s) : 0 \leq s \leq t\right) \text{ and } G_t^{H,i} = \overline{span}\left(B_H^i(s) : 0 \leq s \leq t\right),$$

where $\{B_H^o(t)\}_{t\geq 0}$ and $\{B_H^e(t)\}_{t\geq 0}$ are the odd and even part of $\{B_H(t)\}_{t\in\mathbb{R}}$.

For $t \geq 0$, $i \in \{o, e\}$ and $H + H' = 1$ define

(38) $$M_H^i(t) = \mathbb{E}\left(B_{H'}^i(t) \mid \mathcal{F}_t^{H,i}\right).$$

**Theorem 9.** $\{M_H^o(t)\}_{t \geq 0}$ and $\{M_H^e(t)\}_{t \geq 0}$ are $H$-self-similar Gaussian martingales adapted to the filtration $\{\mathcal{F}_t^{H,o}\}_{t \geq 0}$ and $\{\mathcal{F}_t^{H,e}\}_{t \geq 0}$ respectively.

*Proof.* It is enough to verify the statement for $M_H^o$ only, because the verification for $M_H^e$ is similar. By construction $M_H^o$ is a Gaussian process. It follows from (31) that for $s \leq t$,

$$B_{H'}^o(t) - B_{H'}^o(s) \perp G_s^{H,o},$$

which implies

$$\mathbb{E}\left(M_H^o(t) \mid \mathcal{F}_s^{H,o}\right) = \mathbb{E}\left(\mathbb{E}\left(B_{H'}^o(t) \mid \mathcal{F}_t^{H,o}\right) \mid \mathcal{F}_s^{H,o}\right) = \mathbb{E}\left(B_{H'}^o(t) \mid \mathcal{F}_s^{H,o}\right)$$
$$= \mathbb{E}\left(B_{H'}^o(t) - B_{H'}^o(s) \mid \mathcal{F}_s^{H,o}\right) + \mathbb{E}\left(B_{H'}^o(s) \mid \mathcal{F}_s^{H,o}\right)$$
$$= 0 + M_H^o(s) = M_H^o(s).$$

By $H$-self-similarity property of the odd part $\{B_H^o(t)\}_{t \in [0,\infty), H \in (0,1)}$ of the fBf $\{B_H(t)\}_{t \in \mathbb{R}, H \in (0,1)}$ the field $\{Z_H^o(t)\}_{t \in [0,\infty), H \in (0,1)}$ defined by

$$Z_H^o(t) = a^{-H} B_H^o(at)$$

is an odd fBf and, therefore for $H + H' = 1$,

$$\left\{\mathbb{E}\left(Z_{H'}^o(t) \mid \sigma\left(Z_H^o(r) : 0 \leq r \leq t\right)\right)\right\}_{t \geq 0} \stackrel{f.d.d.}{=} \{M_H^o(t)\}_{t \geq 0}.$$

Furthermore

$$\sigma\left(Z_H^o(r) : 0 \leq r \leq t\right) = \sigma\left(a^{-H} B_H^o(ar) : 0 \leq r \leq t\right) = \sigma\left(B_H^o(ar) : 0 \leq r \leq t\right)$$
$$= \sigma\left(B_H^o(s) : 0 \leq s \leq at\right) = \mathcal{F}_{at}^{H,o}.$$

Therefore

$$\mathbb{E}\left(Z_{H'}^o(t) \mid \sigma\left(Z_H^o(r) : 0 \leq r \leq t\right)\right) = \mathbb{E}\left(a^{-H} B_{H'}^o(at) \mid \mathcal{F}_{at}^{H,o}\right),$$

which concludes the proof. □

So far we have shown that $M_H^o$ and $M_H^e$ are $H$-self-similar Gaussian martingales. By construction, for $i \in \{o, e\}$, $M_H^i(t)$ is an element of $G_t^{H,i}$, and therefore it may be possible to express it as a stochastic integral, up to time $t$, of $B_H^i$. In the case $H = \frac{1}{2}$ this is trivial, since then $H' = \frac{1}{2}$ and therefore $B_H^i = B_{H'}^i$ is a constant multiple of Brownian motion and $M_H^i(t) = B_H^i(t)$. The case $H \neq \frac{1}{2}$ has been solved in our paper [4]. We state the result below without proof. The supporting materials are too long for the present paper. It should also be mentioned that the natural filtration of the martingale $M_H^i$ coincides with the natural filtration of the process $B_H^i$ [4].

**Theorem 10.** Let $H \in (0,1) \setminus \{\frac{1}{2}\}$. If $t \geq 0$ then

$$M_H^o(t) = \frac{\sqrt{\pi} \alpha_H}{\Gamma(1-H)} \int_0^t (t^2 - s^2)_+^{\frac{1}{2}-H} dB_H^o(s)$$

*and*

$$M_H^e(t) = -\frac{\alpha_H}{\Gamma\left(\frac{3}{2} - H\right)} \int_0^t \frac{d}{ds}\left(\int_s^t \left(x^2 - s^2\right)^{\frac{1}{2} - H} dx\right) dB_H^e(s),$$

*where*

$$\alpha_H = \frac{2^{2H-1}\sqrt{\Gamma(3 - 2H)}\sin(\pi H)}{\Gamma\left(\frac{3}{2} - H\right)\sqrt{\Gamma(2H + 1)}}.$$

## 6. Remarks on multifractional Brownian motions

Lévy Véhel and Peltier [23] and Benassi et al. [2] have introduced independently, multifractional Brownian motion. In this section we will clarify the relationship between multifractional Brownian motion and the fractional Brownian fields introduced in Sections 3 and 4. Additionally we will point out an error in the covariance of multifractional Brownian motion obtained from the nonanticipating moving average representation of fBm which shows that in fact the processes of Lévy Véhel [23] and Peltier and Benassi et al. [2] are not the same, as it has been claimed in Cohen [3].

Let $\{W_s\}_{t \in \mathbb{R}}$ be a Brownian motion. For $t \geq 0$ Lévy Véhel and Peltier [23] called

$$(39) \quad X_t = \frac{1}{\left(H(t) + \frac{1}{2}\right)} \int_{\mathbb{R}} \left((t-s)_+^{H(t) - \frac{1}{2}} - (-s)_+^{H(t) - \frac{1}{2}}\right) dW_s,$$

where $H : [0, \infty) \to (0, 1)$ is a deterministic Hölder function with exponent $\beta > 0$, a multifractional Brownian motion. This process is introduced as a generalization of fBm that has different regularity at each $t$, more precisely, if $0 < H(t) < \min(1, \beta)$ then at each $t_0$ the multifractional Brownian motion has Hölder exponent $H(t_0)$ with probability 1. It is clear that if $H(t) \equiv H$ for some $0 < H < 1$, then $\{X_t\}_{t \geq 0}$ is a (nonstandard) $H$-fBm.

Benassi et al. [2] have introduced the process

$$(40) \quad Y_t = \int_{\mathbb{R}} \frac{e^{it\xi} - 1}{|\xi|^{\frac{1}{2} + H(t)}} d\widetilde{\overline{W}}_\xi,$$

where in "some sense" the random measure $d\widetilde{\overline{W}}$ is the Fourier transform of $dW$ and for $g, h \in L^2(\mathbb{R})$ it satisfies

$$\mathbb{E}\left(\int_{-\infty}^{\infty} g(\xi) d\widetilde{\overline{W}}_\xi \overline{\int_{-\infty}^{\infty} h(\xi) d\widetilde{\overline{W}}_\xi}\right) = \int_{-\infty}^{\infty} g(\xi) \overline{h(\xi)} d\xi$$

(see section 7.2.2 in Samorodnitsky and Taqqu [21]). If $H(t) \equiv H$, for some $0 < H < 1$, the process $\{Y_t\}_{t \geq 0}$ is an $H$-fBm, because the right-hand-side of (40) reduces to the well-known harmonizable representation of fBm. The result concerning the Hölder exponent for $\{X_t\}_{t \geq 0}$ holds for the process $\{Y_t\}_{t \geq 0}$ too. Although the process $\{Y_t\}_{t \geq 0}$ is called multifractional Brownian motion, we will refer to it as a harmonizable multifractional Brownian motion to emphasize the differences between the two processes.

In [3] Cohen states that both multifractional Brownian motions $\{X_t\}_{t \geq 0}$ and $\{Y_t\}_{t \geq 0}$ if normalized appropriately are versions of the same process. More precisely the following is stated in [3] (as Theorem 1):

The harmonizable representation of the multifractional Brownian motion:

$$\text{(41)} \qquad \int_{\mathbb{R}} \frac{e^{it\xi} - 1}{|\xi|^{\frac{1}{2}+H(t)}} d\overline{\widehat{W}}_\xi,$$

is almost surely equal up to a multiplicative deterministic function to the well balanced moving average

$$\int_{\mathbb{R}} \left(|t-s|^{H(t)-\frac{1}{2}} - |s|^{H(t)-\frac{1}{2}}\right) dW_s.$$

When $H(t) = \frac{1}{2}$, $\left(|t-s|^{H(t)-\frac{1}{2}} - |s|^{H(t)-\frac{1}{2}}\right)$ is ambiguous, hence the conventional meaning

$$|t-s|^0 - |s|^0 = \log\left(\frac{1}{|t-s|}\right) - \log\left(\frac{1}{|s|}\right)$$

is to be used. Conversely, one can show that the non anticipating moving average

$$\text{(42)} \qquad \int_{\mathbb{R}} \left((t-s)_+^{H(t)-\frac{1}{2}} - (s)_+^{H(t)-\frac{1}{2}}\right) dW_s$$

is equal up to a multiplicative deterministic function to the harmonizable representation

$$\int_{\mathbb{R}} \frac{e^{it\xi} - 1}{i\xi |\xi|^{H(t)-\frac{1}{2}}} d\overline{\widehat{W}}_\xi.$$

Hence the mfBm given by the non anticipating moving average (42) has the same law as the mfBm given by the harmonizable representation (41) up to a multiplicative deterministic function.

The arguments used in the proof of the of Theorem 1 in Cohen [3] are based on the fact that the Fourier transform of $x \to |t-x|^{H(t)-\frac{1}{2}} - |x|^{H(t)-\frac{1}{2}}$ for $H(t) \neq \frac{1}{2}$ and of $\log(\frac{1}{|t-x|}) - \log(\frac{1}{|x|})$ for $H(t) = \frac{1}{2}$ equal, up to a multiplicative constant, to $\xi \to \frac{e^{it\xi}-1}{|\xi|^{\frac{1}{2}+H(t)}}$; and an incorrect statement that the Fourier transform of $x \to (t-x)_+^{H(t)-\frac{1}{2}} - (-x)_+^{H(t)-\frac{1}{2}}$ equals, up to a multiplicative constant, to $\xi \to \frac{e^{it\xi}-1}{i\xi|\xi|^{H(t)-\frac{1}{2}}}$. The equations (16) and (18) are the correct expression for that Fourier transform. Consequently the last two statements of the above Theorem 1 in Cohen are incorrect. In order to see why, consider two multifractional Brownian motions

$$X_t = \frac{1}{c_{H(t)}} \int_{-\infty}^{\infty} (t-x)_+^{H(t)-\frac{1}{2}} - (-x)_+^{H(t)-\frac{1}{2}} \, dB_x$$

and

$$Y_t = \begin{cases} \frac{1}{d_{H(t)}} \int_{-\infty}^{\infty} |t-x|^{H(t)-\frac{1}{2}} - |x|^{H(t)-\frac{1}{2}} \, dW_x & \text{for } t \, H(t) \neq \frac{1}{2} \\ \frac{1}{d_{H(t)}} \int_{-\infty}^{\infty} \log\left(\frac{1}{|t-x|}\right) - \log\left(\frac{1}{|x|}\right) dW_x & H(t) = \frac{1}{2} \end{cases},$$

where $t \geq 0$, $c_{H(t)}$ and $d_{H(t)}$ are defined by (6), (8) respectively, $d_{\frac{1}{2}} = \pi$, and $\{B_t\}_{t\in\mathbb{R}}$ and $\{W_t\}_{t\in\mathbb{R}}$ are Brownian motions. According to the last statement of Theorem 1 in Cohen [3] there is a deterministic function $f_t$ such that the processes $\{X_t\}_{t\geq 0}$ and $\{f_t Y_t\}_{t\geq 0}$ have the same law. The chosen normalization assures that $\mathbb{E}(X_t^2) = \mathbb{E}(Y_t^2)$ for all $t \geq 0$, implying that $|f_t| = 1$ for all $t$ such that $t > 0$. It follows that $|\mathbb{E}(X_t X_s)| = |\mathbb{E}(Y_t Y_s)|$ for all $s,t$. It is clear that the process $\{X_t\}_{t\geq 0}$ can be obtained from a dependent fractional Brownian field $\{B_H(t)\}_{t\geq 0, H\in(0,1)}$ as $\{B_{H(t)}(t)\}_{t\geq 0}$. Similarly, the Gaussian field $\{W_H(t)\}_{t\geq 0, H\in(0,1)}$ defined for $H\neq \frac{1}{2}$ by equation (7) and for $H=\frac{1}{2}$ by (9) gives $\{Y_t\}_{t\geq 0}$ via $\{W_{H(t)}(t)\}_{t\geq 0}$. Since the last statement of Theorem 1 in [3] is supposed to hold for every Hölder function $t\to H(t)$, that statement holds if and only if the Gaussian fields $\{B_H(t)\}_{t\in[0,\infty), H\in(0,1)}$ and $\{W_H(t)\}_{t\in[0,\infty), H\in(0,1)}$ have the same covariance in absolute value. The covariance of $\{B_H(t)\}_{t\in[0,\infty), H\in(0,1)}$ is given by Theorem 2. Proposition 11 right below gives the covariance of $\{W_H(t)\}_{t\in[0,\infty), H\in(0,1)}$. The proposition shows that if $H\neq \frac{1}{2}$ and $H+H'=1$ the covariance $\mathbb{E}W_H(t)W_{H'}(s)$ is a multiple of $s\wedge t$. From Theorem 2 we see that this is not the case for $\mathbb{E}B_H(t)B_{H'}(s)$.

**Proposition 11.** *The covariance of the Gaussian field $\{W_H(t)\}_{t\in\mathbb{R}, H\in(0,1)}$, defined for $H\neq \frac{1}{2}$ by equation (7) and for $H=\frac{1}{2}$ by (9) is*

$$\mathbb{E}W_H(t)W_{H'}(s) = \frac{k_H k_{H'}}{d_H d_{H'}} \frac{d^2_{\frac{H+H'}{2}}}{k^2_{\frac{H+H'}{2}}} \left( \frac{|t|^{H+H'} + |s|^{H+H'} - |t-s|^{H+H'}}{2} \right),$$

*where $d_H$ is defined by (8) and*

(43) $$k_H = -2\Gamma\left(H+\frac{1}{2}\right)\sin\left(\left(H-\frac{1}{2}\right)\frac{\pi}{2}\right) \text{ for } H\neq \frac{1}{2} \text{ and } k_{\frac{1}{2}}=\pi.$$

Before proving Proposition 11, a few technical results are needed.

**Lemma 12.** *Let $f_{t,\frac{1}{2}}(x) = \log\frac{1}{|t-x|} - \log\frac{1}{|x|}$ for $x,t\in\mathbb{R}$. Then*

$$\widehat{f_{t,\frac{1}{2}}}(\xi) = \pi\frac{e^{it\xi}-1}{\xi}$$

*Proof.* Suppose that $\xi \neq 0$. Since $f_{t,\frac{1}{2}} \in L^1(\mathbb{R}) \cap L^2(\mathbb{R})$ the Fourier transform of $f_{t,\frac{1}{2}}$ can be computed as

$$\widehat{f_{t,\frac{1}{2}}}(\xi) = \int_{-\infty}^{\infty} \left(\log\frac{1}{|t-x|} - \log\frac{1}{|x|}\right) e^{ix\xi} dx = \int_{-\infty}^{\infty} \log\left(\frac{|x|}{|x-t|}\right) e^{ix\xi} dx.$$

Substituting $u = x - \frac{t}{2}$ yields

$$\widehat{f_{t,\frac{1}{2}}}(\xi) = e^{i\frac{t}{2}\xi} \int_{-\infty}^{\infty} \log\left(\frac{|u+\frac{t}{2}|}{|u-\frac{t}{2}|}\right) e^{iu\xi} du,$$

and after the substitution $u=-v$ on $(-\infty, 0)$ to

$$\widehat{f_{t,\frac{1}{2}}}(\xi) = e^{i\frac{t}{2}\xi} 2i \int_0^{\infty} \log\left(\frac{|u+\frac{t}{2}|}{|u-\frac{t}{2}|}\right) \sin(u\xi) du.$$

Formula 4.382 (1) in [7] reads:

$$\int_0^{\infty} \log\left(\frac{|u+a|}{|u-a|}\right) \sin(bx) dx = \frac{\pi}{b}\sin(ab) \text{ for } a,b>0.$$

If $t>0$ set $a=\frac{t}{2}$ and $b=|\xi|$ and use that $\sin(z\xi) = sgn(\xi)\sin(z|\xi|)$ for $z\geq 0$ to get

$$\widehat{f_{t,\frac{1}{2}}}(\xi) = e^{i\frac{t}{2}\xi} 2i \frac{\pi}{|\xi|} \sin\left(\frac{t\xi}{2}\right)$$

Standard trigonometric identities $\sin(2a) = 2\sin(a)\cos(a)$ and $\sin^2(a) = \frac{1-\cos(2a)}{2}$ complete the proof for $t > 0$. The proof for $t < 0$ is similar. $\square$

**Lemma 13.** Let $f_{t,H}(x) = |t-x|^{H-\frac{1}{2}} - |x|^{H-\frac{1}{2}}$ for $x, t \in \mathbb{R}$ and $H \in (0,1) \setminus \{\frac{1}{2}\}$. Then
$$\widehat{f_{t,H}}(\xi) = -2\Gamma\left(H + \frac{1}{2}\right) \sin\left(\left(H - \frac{1}{2}\right)\frac{\pi}{2}\right) \frac{e^{it\xi} - 1}{|\xi|^{H+\frac{1}{2}}}.$$

*Proof.* Rewriting $f_{t,H}$ as
$$f_{t,H}(x) = (t-x)_+^{H-\frac{1}{2}} - (x)_+^{H-\frac{1}{2}} + (t-x)_-^{H-\frac{1}{2}} - (x)_-^{H-\frac{1}{2}},$$

and applying equations (16) and (17) it follows that
$$\widehat{f_{t,H}}(\xi) = \Gamma\left(H + \frac{1}{2}\right) \frac{e^{it\xi} - 1}{i\xi} (i\xi)^{-(H-\frac{1}{2})} - \Gamma\left(H + \frac{1}{2}\right) \frac{e^{it\xi} - 1}{i\xi} (-i\xi)^{-(H-\frac{1}{2})}.$$

and from (18) that
$$(i\xi)^{-(H-\frac{1}{2})} - (-i\xi)^{-(H-\frac{1}{2})} = 2i|\xi|^{-(H-\frac{1}{2})} sgn(-\xi) \sin\left(\left(H - \frac{1}{2}\right)\frac{\pi}{2}\right).$$

Therefore
$$\widehat{f_{t,H}}(\xi) = 2\Gamma\left(H + \frac{1}{2}\right) \frac{e^{it\xi} - 1}{i\xi} i|\xi|^{-(H-\frac{1}{2})} sgn(-\xi) \sin\left(\left(H - \frac{1}{2}\right)\frac{\pi}{2}\right)$$
$$= -2\Gamma\left(H + \frac{1}{2}\right) \sin\left(\left(H - \frac{1}{2}\right)\frac{\pi}{2}\right) \frac{e^{it\xi} - 1}{|\xi|^{H+\frac{1}{2}}}.$$
$\square$

We are now ready to prove Proposition 11. The idea is to use the fact that up to a multiplicative constant $\widehat{f_{t,H}}(\xi)\overline{\widehat{f_{s,H'}}(\xi)}$ equals $\widehat{f_{t,\frac{H+H'}{2}}}(\xi)\overline{\widehat{f_{s,\frac{H+H'}{2}}}(\xi)}$ and that up to a multiplicative constant the integral over $\mathbb{R}$ of the later is the covariance of an $\frac{H+H'}{2}$-fBm. This argument is used in Ayache et al. [1] to compute the covariance of a multifractional Brownian motion given by (40). In Ayache et al. [1] it is also erroneously claimed that the covariance of the multifractional Brownian motion given by (39) is the same, if properly normalized, as the covariance of the harmonizable multifractional Brownian motion given by (40). Their proof is based on the last statement of Theorem 1 in Cohen [3]. In section 3 of Lim and Muniandy [11], the authors give another incorrect argument about the equivalence (up to a deterministic multiplicative function) between the harmonizable multifractional Brownian motion (40) and the nonanticipative multifractional Brownian motion (39).

*Proof of Proposition 11.* By the Plancherel identity
(44)
$$\mathbb{E}W_H(t)W_{H'}(s) = \frac{1}{d_H d_{H'}} \int_{-\infty}^{\infty} f_{t,H}(x) f_{s,H'}(x) dx$$
$$= \frac{1}{2\pi} \frac{1}{d_H d_{H'}} \int_{-\infty}^{\infty} \widehat{f_{t,H}}(\xi)\overline{\widehat{f_{s,H'}}(\xi)} d\xi,$$

where $f_{t,H}$, when $H = \frac{1}{2}$, is the same as in Lemma 12, and, when $H \neq \frac{1}{2}$, as in and Lemma 13. From Lemma 12 and Lemma 13 it follows

$$\int_{-\infty}^{\infty} \widehat{f_{t,H}}(\xi) \overline{\widehat{f_{s,H'}}(\xi)} d\xi = k_H k_{H'} \int \frac{(e^{it\xi}-1)\overline{(e^{is\xi}-1)}}{|\xi|^{H+\frac{1}{2}} |\xi|^{H'+\frac{1}{2}}} d\xi$$

Let $H_0 = \frac{H+H'}{2}$. Observe that $H_0 \in (0,1)$ and so

$$\int_{-\infty}^{\infty} \widehat{f_{t,H}}(\xi) \overline{\widehat{f_{s,H'}}(\xi)} d\xi = k_H k_{H'} \frac{d_{H_0}^2}{k_{H_0}^2} \frac{k_{H_0}^2}{d_{H_0}^2} \int \frac{(e^{it\xi}-1)\overline{(e^{is\xi}-1)}}{|\xi|^{H_0+\frac{1}{2}} |\xi|^{H_0+\frac{1}{2}}} d\xi$$

$$= k_H k_{H'} \frac{d_{H_0}^2}{k_{H_0}^2} \frac{2\pi}{d_{H_0}^2} \int_{-\infty}^{\infty} f_{t,H_0}(x) f_{s,H_0}(x) dx,$$

where the last equality follows from Lemma 12, Lemma 13 and Plancherel identity. Hence

$$\int_{-\infty}^{\infty} \widehat{f_{t,H}}(\xi) \overline{\widehat{f_{s,H'}}(\xi)} d\xi = k_H k_{H'} \frac{d_{H_0}^2}{k_{H_0}^2} 2\pi \mathbb{E} W_{H_0}(t) W_{H_0}(s)$$

$$= k_H k_{H'} \frac{d_{H_0}^2}{k_{H_0}^2} 2\pi \left( \frac{|t|^{2H_0} + |s|^{2H_0} - |t-s|^{2H_0}}{2} \right).$$

Hence (44) becomes

$$\mathbb{E} W_H(t) W_{H'}(s) = \frac{k_H k_{H'}}{d_H d_{H'}} \frac{d_{H_0}^2}{k_{H_0}^2} \left( \frac{|t|^{2H_0} + |s|^{2H_0} - |t-s|^{2H_0}}{2} \right).$$

This finishes the proof. $\square$